\documentclass[11pt]{amsart}
\usepackage{amsfonts}
\usepackage{amssymb,amsfonts,amsmath,bm}

\begin{document}
\theoremstyle{plain}
\newtheorem{thm}{Theorem}[section]
\newtheorem{theorem}[thm]{Theorem}
\newtheorem{lemma}[thm]{Lemma}
\newtheorem{corollary}[thm]{Corollary}
\newtheorem{corollary*}[thm]{Corollary*}
\newtheorem{proposition}[thm]{Proposition}
\newtheorem{proposition*}[thm]{}
\newtheorem{conjecture}[thm]{Conjecture}
\theoremstyle{definition}
\newtheorem{construction}[thm]{Construction}
\newtheorem{notations}[thm]{Notations}
\newtheorem{question}[thm]{Question}
\newtheorem{problem}[thm]{Problem}
\newtheorem{remark}[thm]{Remark}
\newtheorem{remarks}[thm]{Remarks}
\newtheorem{definition}[thm]{Definition}
\newtheorem{claim}[thm]{Claim}
\newtheorem{assumption}[thm]{Assumption}
\newtheorem{assumptions}[thm]{Assumptions}
\newtheorem{properties}[thm]{Properties}
\newtheorem{example}[thm]{Example}
\newtheorem{comments}[thm]{Comments}
\newtheorem{blank}[thm]{}
\newtheorem{observation}[thm]{Observation}
\newtheorem{defn-thm}[thm]{Definition-Theorem}

\newcommand{\sM}{{\mathcal M}}

\def\arcsinh{\operatorname{arcsinh}}
\def\rank{\operatorname{rank}}
\def\Res{\operatorname{Res}}


\title[Computing top intersections in $\mathcal M_g$]{Computing top intersections in the \\tautological ring of $\mathcal M_g$}
        \author{Kefeng Liu}
        \address{Center of Mathematical Sciences, Zhejiang University, Hangzhou, Zhejiang 310027, China;
                Department of Mathematics,University of California at Los Angeles,
                Los Angeles, CA 90095-1555, USA}
        \email{liu@math.ucla.edu, liu@cms.zju.edu.cn}
        \author{Hao Xu}
        \address{Center of Mathematical Sciences, Zhejiang University, Hangzhou, Zhejiang 310027, China;
               Department of Mathematics, Harvard University, Cambridge, MA 02138, USA}
        \email{haoxu@math.harvard.edu}
        

        \begin{abstract}
        We derive effective recursion formulae
        of top intersections in the tautological ring $R^*(\mathcal
        M_g)$ of the moduli space of curves of genus $g\geq 2$. As an application, we prove
        a convolution-type tautological relation in $R^{g-2}(\mathcal
        M_g)$.
        \end{abstract}
    \maketitle

\section{Introduction}

We denote by $\overline{\sM}_{g,n}$ the moduli space of stable
$n$-pointed genus $g$ complex algebraic curves. We have the morphism
that forgets the last marked point,
$$
\pi_{n+1}: \overline{\sM}_{g,n+1}\longrightarrow
\overline{\sM}_{g,n}.
$$
Denote by $\sigma_1,\dots,\sigma_n$ the canonical sections of $\pi$,
and by $D_1,\dots,D_n$ the corresponding divisors in
$\overline{\sM}_{g,n+1}$. Let $\omega_{\pi}$ be the relative
dualizing sheaf. We have the following tautological classes on
moduli spaces of curves.
\begin{align*}
\psi_i&=c_1(\sigma_i^*(\omega_{\pi})),\\
\kappa_i&=\pi_*\left(c_1\left(\omega_{\pi}\left(\sum D_i\right)\right)^{i+1}\right),\\
\lambda_l&=c_l(\pi_*(\omega_{\pi})),\quad 1\leq l\leq g.
\end{align*}
The definition of $\kappa$ classes on $\overline{\sM}_{g,n}$ is due
to Arbarello-Cornalba \cite{AC}, generalizing Mumford-Morita-Miller
classes
$$
\kappa_i=\pi_*\left(c_1\left(\omega_{\pi}\right)^{i+1}\right)\in
\mathcal A^i(\mathcal M_g),
$$
where $\mathcal A^*(\mathcal M_g)$ is the rational Chow ring of
$\mathcal M_g$.

The tautological ring $R^*(\mathcal M_g)$ is defined to be the
subalgebra of $\mathcal A^*(\mathcal M_g)$ generated by the
tautological classes $\kappa_i$.

We use Witten's notation to denote intersection numbers:
$$\langle\tau_{d_1}\cdots\tau_{d_n}\kappa_{a_1}\cdots\kappa_{a_m}\mid\lambda_{1}^{k_{1}}
    \cdots\lambda_{g}^{k_{g}}\rangle:= \int_{\overline{\mathcal{M}}_{g,n}}\psi_{1}^{d_{1}}\cdots\psi_{n}^{d_{n}}\kappa_{a_1}\cdots\kappa_{a_m}\lambda_{1}^{k_{1}}
    \cdots\lambda_{g}^{k_{g}}.$$
These are rational numbers and are called the Hodge integrals, which
can be computed by Faber's algorithm \cite{Fa2} based on Mumford's
formula for Chern characters of Hodge bundles and the celebrated
Witten-Kontsevich theorem \cite{Wi, Ko}.

\subsection{Faber's conjecture} Around 1993, Faber \cite{Fa} proposed a series of
remarkable conjectures about the structure of $R^*(\sM_g)$:

\begin{enumerate}
\item[i)] For $0 \leq k \leq g-2$, the natural
product $$R^k(\mathcal M_g) \times R^{g-2-k}(\mathcal M_g)
\rightarrow R^{g-2}(\mathcal M_g) \cong \mathbb Q$$ is a perfect
pairing.

\item[ii)] The $[g/3]$ classes $\kappa_1,\dots,\kappa_{[g/3]}$
generate the ring $R^*(\sM_g)$, with no relations in degrees $\leq
[g/3]$.

\item[iii)] Let $\sum_{j=1}^n d_j=g-2$ and $d_j\geq0$. Then
\begin{equation}
\sum_{\sigma\in
S_n}\kappa_\sigma=\frac{(2g-3+n)!}{(2g-2)!!\prod_{j=1}^{n}(2d_j+1)!!}\kappa_{g-2},
\end{equation}
where $\kappa_\sigma$ is defined as follows: write the permutation
$\sigma$ as a product of $\nu(\sigma)$ disjoint cycles
$\sigma=\beta_1\cdots\beta_{\nu(\sigma)}$, where we think of the
symmetric group $S_n$ as acting on the $n$-tuple $(d_1,\dots ,d_n)$.
Denote by $|\beta|$ the sum of the elements of a cycle $\beta$. Then
$ \kappa_\sigma=\kappa_{|\beta_1|}\kappa_{|\beta_2|}\dots
\kappa_{|\beta_{\nu(\sigma)}|}$.

\end{enumerate}

It is a theorem of Looijenga \cite{Lo} that
$$\dim R^k(\mathcal M_g)=0,\quad k>g-2,$$
$$\dim R^{g-2}(\mathcal M_g)\leq1.$$ Faber proved that
actually $\dim R^{g-2}(\mathcal M_g)=1$.

Part (i) of Faber's conjecture is also called Faber's perfect
pairing conjecture, which is still open. Faber has verified
$g\leq23$.

Part (ii) has been proved independently by Morita \cite{Mo} and
Ionel \cite{Io} by very different methods. As pointed out by Faber
\cite{Fa}, Harer's stability result implies that there is no
relation in degrees $\leq [g/3]$.

Part (iii) is known as the Faber intersection number conjecture and
is equivalent to
\begin{equation}\label{eqfaber}
\langle\tau_{d_1+1}\cdots\tau_{d_n+1}\mid\lambda_g\lambda_{g-1}\rangle_g
=\frac{(2g-3+n)!}{(2g-2)!!\prod_{j=1}^{n}(2d_j+1)!!}\langle\kappa_{g-2}\mid\lambda_g\lambda_{g-1}\rangle_g.
\end{equation}
A short and direct proof of the Faber intersection number conjecture
can be found in \cite{LX2}. Also see \cite{GP, GJV} and the most
recently \cite{BS} for different approaches to the problem.

In fact, Faber \cite{Fa} further proposed that the tautological ring
$R^*(\mathcal M_g)$ behaves like the algebraic cohomology ring of a
nonsingular projective variety of dimension $g-2$, i.e. it satisfies
the Hard Lefschetz and Hodge Positivity properties with respect to
$\kappa_1$.

Faber's intersection number conjecture determines the top
intersections in $R^{g-2}(\mathcal M_g)$. If we assume Faber's
perfect pairing conjecture, then the ring structure of $\mathcal
R^{*}(\mathcal M_g)$ is also determined.

A central theme in Faber's conjecture is to explicitly describe
relations in the tautological rings. However, it is a highly
nontrivial task to identify tautological relations in $R^*(\sM_g)$
when $g$ becomes larger. In this paper, we will consider only
tautological relations of top degree in $R^{g-2}(\mathcal M_g)$.
Already known examples include Faber-Zagier's formula \cite{Fa}
\begin{equation}\label{eqfaza2}
\kappa_1^{g-2}=\frac{1}{g-1}2^{2g-5}((g-2)!)^2\kappa_{g-2}.
\end{equation}
and
Pandharipande's formula \cite{Pa}
\begin{equation}\label{eqpand}
\sum_{i=0}^{g-2}(-1)^i\lambda_i
\kappa_{g-2-i}=\frac{2^{g-1}}{g!}\kappa_{g-2}.
\end{equation}
in $R^{g-2}(\mathcal M_g), g\geq2$.

Now we describe the main results of this paper.
We will prove two effective recursive formulae of different flavors
for computing top intersections in $R^{g-2}(\mathcal M_g)$ (see
Theorems \ref{const5} and \ref{recur1}).  For example, Theorem \ref{recur1} 
can be equivalently stated as the following:
\begin{theorem} \label{recur2}
Let $g\geq 3$ and $|\bold m|=g-2$. Then the following relation
\begin{equation} \label{eqconst12}
\kappa(\bold m)=\frac{1}{(||\bold m||-1)}\sum_{\substack{\bold L+\bold{L'}=\bold m\\||\bold{L}||\geq2}}A_{g,\bold{L}}\binom{\bold m}{\bold L}
\kappa(\bold{L'}+\bm\delta_{|\bold{L}|})
\end{equation}
holds in $R^{g-2}(\mathcal M_g)$, where $A_{g,\bold L}=\bold L! D_{g,\bold L}$ are some
explicitly known constants.
\end{theorem}

Note that in the right-hand side of equation \eqref{eqconst12}, $||\bold{L'}+\bm\delta_{|\bold{L}|}||<||\bold m||$, so equation \eqref{eqconst12} is
indeed an effective recursion relation. Our strategy of proof is to exploit
the method used in \cite{LX1}. From an algorithmic point of view,
recursion formulae are often more effective than closed formulae,
since previously computed values can be reused in recursive computations.

Our recursion formulae will be used to compute Faber's intersection
matrix, whose rank is equal to the dimension of $R^{*}(\mathcal M_g)$ by Faber's perfect pairing conjecture. On the other hand, the cohomological dimension
$H^*(\mathcal M_g)$ is also an outstanding open problem \cite{AC2}.

One objective of this paper is to gain a better understanding of Faber-Zagier's formula.
In Section \ref{Bern}, we prove an interesting Bernoulli number identity equivalent to
Faber-Zagier's formula.

In the final section, we prove the following convolution-type tautological
relation:
\begin{theorem} \label{tauto} Let $g\geq3$. We have the following relation in $R^{g-2}(\mathcal
M_g)$,
\begin{equation}
\sum_{i=0}^{g-2}D_{g,g-2-i}\frac{\kappa_1^i \kappa_{g-2-i}}{i!}=0,
\end{equation}
where $D_{g,k}$ are given by
\begin{equation}\label{eqconst10}
D_{g,k}=\frac{3}{2(g-2)}\cdot\frac{1}{k!}+\frac{2g-1}{2(g-2)}\sum^k_{j=0}\frac{(2j+1)(-1)^{j+1}2^{2j}B_{2j}}{j!(k-j)!}.
\end{equation}
where $B_{2j}$ is the $2j$-th Bernoulli number.
We have $D_{g,0}=-1,\ D_{g,1}=\frac{g+1}{g-2},\
D_{g,2}=\frac{17g-4}{6(g-2)}$.

\end{theorem}

We don't know whether the expression of $D_{g,k}$ can be simplified.

\

\noindent{\bf Acknowledgements.}
We would like to thank Professor Jian Zhou for helpful discussions and for kindly providing
an elegant proof of Lemma \ref{zhou}.

\vskip 30pt
\section{The Faber intersection matrix} \label{sectiontr}
First we fix notation. Consider the semigroup $N^\infty$ of
sequences ${\bold m}=(m(1),m(2),\dots)$ where $m(i)$ are nonnegative
integers and $m(i)=0$ for sufficiently large $i$. We sometimes also
use $(1^{m(1)}2^{m(2)}\dots)$ to denote $\bold m$.

Let $\bold m, \bold{a_1,\dots,a_n} \in N^\infty$, $\bold
m=\sum_{i=1}^n \bold{a_i}$.
$$|\bold m|:=\sum_{i\geq 1}i m(i)\quad ||\bold m||:=\sum_{i\geq1}m(i) \quad \binom{\bold m}{\bold{a_1,\dots,a_n}}:=\prod_{i\geq1}\binom{m(i)}{a_1(i),\dots,a_n(i)}.$$

Let $\bold m\in N^\infty$, we denote a formal monomial of $\kappa$
classes by
$$\kappa(\bold m):=\prod_{i\geq1}\kappa_i^{m(i)}.$$

If $|\bold m|=g-2$, then from Faber's intersection number identity
\eqref{eqfaber} and the formula expressing $\psi$ classes by
$\kappa$ classes \cite{KMZ}, we have the following relation in
$R^{g-2}(\mathcal M_g)$,
\begin{equation}\label{eqkmz}
\kappa(\bold m)={\rm Fab}_g(\bold m)\kappa_{g-2},
\end{equation}
where the proportional constant ${\rm Fab}_g(\bold m)$ is given by
$${\rm Fab}_g(\bold m)=\sum_{r=1}^{||\bold m||}\frac{(-1)^{||\bold
m||-r}}{r!}\sum_{\substack {\bold
m=\bold{m_1}+\cdots+\bold{m_r}\\\bold{m_i}\neq\bold 0}}\binom{\bold
m}{\bold
{m_1,\dots,m_r}}\frac{(2g-3+r)!}{(2g-2)!!\prod_{j=1}^r(2|\bold
m_j|+1)!!}.$$

Let $g\geq 2$ and $0\leq k\leq g-2$. Denote by $p(n)$ the number of
partitions of $n$. Define a matrix $V_g^k$ of size $p(k)\times
p(g-2-k)$ with entries
\begin{equation}\label{eqkmz2}
(V_g^k)_{\bold L, \bold L'}={\rm Fab}_g(\bold L+\bold{L'}),
\end{equation}
where $\bold L, \bold{L'}\in N^\infty$ and $|\bold L|=k$, $|\bold
L'|=g-2-k$.

We call $V_g^k$ the Faber intersection matrix. If Faber's perfect
pairing conjecture is true, then we have
\begin{equation}\label{eqper}
\rank V_g^k=\dim R^k(\sM_g),\qquad 0\leq k\leq g-2.
\end{equation}
 Faber has
verified his conjecture for all $g\leq23$, so the above relation
holds for at least $g\leq23$. Thus, we may get useful information of
$R^*(\sM_g)$ from the Faber intersection matrix.

In the next section, we will develop a recursive method for
computing entries of $V_g^k$. As a result, we have computed $V_g^k$
for all $g\leq36$. These data is used to check a conjectural relationship
between $\rank V_g^k$ and Ramanujan's mock theta function (see Section 7 of \cite{LX3}). 

The recursive algorithm did not reduce
computational complexity in theory. But the recursively reusable
data rendered the computation more efficient.

If we calculate $V_g^k$ by equation \eqref{eqkmz2} on a computer,
the required CPU time is proportional to the number of vector
partitions.

\begin{definition}
For $\bold m\in N^{\infty}$, denote by $P(\bold m)$ the number of
distinct representations of $\bold m$ as an unordered sum in
$N^{\infty}$,
$$\bold m={\bold m}_1+{\bold m}_2+\cdots+{\bold m}_r,$$
where ${\bold m}_i\neq0$. We call such $P(\bold m)$ the vector
partition number.
\end{definition}

We make the convention that $P(\bold 0)=1$. There is a simple
recursion formula for computing $P(\bold m)$, due to Cheema and
Motzkin \cite{CM}.
\begin{lemma} {\bf (Cheema-Motzkin)} For $\bold m=(m_1,m_2,\dots,m_s)$ with
$m_1\neq0$, we have
\begin{equation} \label{eqpar}
m_1P(\bold
m)=\sum_{a_i\geq0}\left(P(m_1-a_1,\dots,m_s-a_s)\sum_{\substack{k>0\\k\mid
\gcd(a_1,\dots,a_s)}}\frac{a_1}{k}\right).
\end{equation}
\end{lemma}
\begin{proof}Let $\bold x=(x_1,x_2,\dots)$ be a sequence of formal variables and
$\bold x^{\bold m}=x_1^{m(1)} x_2^{m(2)}\dots$ for $\bold m\in
N^\infty$. We have
$$\sum_{\bold m\in N^{\infty}}P(\bold m)\bold x^{\bold m}=\prod_{\substack{\bold a\in N^\infty\\\bold a\neq\bold 0}}(1-{\bold x}^{\bold a})^{-1}.$$
Hence
$$\log \sum_{\bold m\in N^{\infty}}P(\bold m)\bold x^{\bold m}=
-\sum_{\substack{\bold a\in N^\infty\\\bold a\neq\bold 0}}
\log(1-\bold x^{\bold a}) =\sum_{\substack{\bold a\in
N^\infty\\\bold a\neq\bold 0}} \sum_{k=1}^\infty\frac{\bold
x^{k\bold a}}{k}$$ Differentiating with respect to $x_1$, we get
$$\sum_{\bold m\in N^{\infty}}m_1 P(\bold m)\bold x^{\bold m}=\sum_{\bold m\in N^{\infty}}P(\bold m)\bold x^{\bold m}
\sum_{\substack{\bold a\in N^\infty\\\bold a\neq\bold
0}}\sum_{k=1}^\infty a_1\bold x^{k\bold a}.$$ Equating coefficients
of $\bold x^{\bold m}$, we get the desired identity \eqref{eqpar}.
\end{proof}

For a given genus $g$, the complexity of computing $V_g^k$ for all
$0\leq k\leq g-2$ is measured by the following quantity
$$D(g-2):=\sum_{|\bold m|=g-2}P(\bold m),$$

Some values of $D(n)$ are listed as follows.
$$
\begin{array}{c||c|c|c|c|c|c|c|c|c|c}
 n & 0 & 1 &2 &3&4&5&10&20&30&40
\\\hline D(n) & 1&1&3&6&14&27&817&318106&71832114&11668071461
\end{array}
$$

\begin{proposition} $D(n)$ has a simple generating function
$$\sum_{n=0}^\infty D(n)x^n=\prod_{n=1}^\infty\frac{1}{(1-x^n)^{p(n)}}.$$
\end{proposition}
\begin{proof}
We have
$$\sum_{\bold m}P(\bold m)\bold x^{\bold m}=\prod_{\bold m\neq\bold 0} \frac{1}{(1-\bold x^{\bold m})}.$$
Substitute $x_i$ by $x^i$, we get the desired result.
\end{proof}

$D(n)$ can also be regarded as the number of double partitions of
$n$. The following asymptotic formula is proved by Kaneiwa
\cite{Kan}
$$D(n)\thicksim e^{\frac{\pi^2 n}{6\log n}},$$
which gives the computational complexity of the Faber intersection
matrix.

Let $\mathcal I(\mathcal M_g)$ be the ideal of polynomial relations
of $\kappa$ classes in $\mathcal M_g$ and let $\mathcal I^k(\mathcal
M_g)$ be the group of relations in degree $k$. We have
$$\dim
R^k(\mathcal M_g)+\dim\mathcal I^k(\mathcal M_g)=p(k).$$

Let $s\geq 0$ and $g=3k-s$. Faber \cite{Fa} pointed out that when
$k\geq s+2$ (i.e. $2k\leq g-2$), $\dim\mathcal I^k(\mathcal M_g)$
depends only on $s$. Denoting this number by $a(s)$, Faber has
\cite{Fa} computed the first ten values. Our calculation of the rank
of $V_g^k$ for $g\leq36$ extends Faber's table of the function $a$.
$$
\begin{array}{c||c|c|c|c|c|c|c|c|c|c|c|c|c|c|c}
 s & 1 & 2 &3 &4&5&6&7&8&9&10&11&12&13&14&15
\\ \hline a(s) & 1&1&2&3&5&6&10&13&18&24&33&41&56&71&91
\end{array}
$$

In \cite{Fa}, Faber says ``Zagier and I have a favourite guess of
this function $a$, but there are many functions with ten prescribed
values.''

We have also tried to guess the function $a$, but failed. For
example, define
$$
f(s)=\sum_{0\leq r\leq [s/3]}p(s+1-3r)-p(s-3r).
$$
We have $f(s)=a(s),\ s\leq10$, but $f(11)=34\neq a(11)$.

\vskip 30pt
\section{Computing top intersections in $R^{g-2}(\mathcal M_g)$}

Let $d_j\geq 1$ and $\sum_j d_j+|\bold m|=g-2+n$. We define the
following quantities
\begin{equation} \label{eqconst7}
F_{g,n}(\bold m):=\frac{(2g-2)!!\prod_{j=1}^n(2d_j-1)!!}{(2g+n-3)!\bold m!}
\cdot\frac{\langle\prod_{j=1}^n\tau_{d_j}\kappa(\bold
m)\mid\lambda_g\lambda_{g-1}\rangle_g}{\langle\kappa_{g-2}\mid\lambda_g\lambda_{g-1}\rangle_g}.
\end{equation}

\begin{proposition}
The above definition of $F_{g,n}(\bold m)$ is independent of $d_j$
and $F_{g,n}(\bold 0)=1$.
\end{proposition}
\begin{proof} From the formula expressing $\psi$ by $\kappa$ in
\cite{KMZ}, we have
\begin{multline*}
\langle\prod_{j=1}^n\tau_{d_j}\kappa(\bold
m)\mid\lambda_g\lambda_{g-1}\rangle_g\\
=\sum_{r=0}^{||\bold
m||}\frac{(-1)^{||\bold m||-r}}{r!}\sum_{\substack {\bold
m=\bold{m_1}+\cdots+\bold{m_r}\\\bold{m_i}\neq\bold 0}}\binom{\bold
m}{\bold
{m_1,\dots,m_r}}\langle\prod_{j=1}^n\tau_{d_j}\prod_{j=1}^r\tau_{|\bold{m_j}|+1}\mid\lambda_g\lambda_{g-1}\rangle_g
\end{multline*}
So the proposition follows directly from Faber's intersection number
identity \eqref{eqfaber}.
\end{proof}

We will also write $F_{g}(\bold m)$ instead of $F_{g,0}(\bold m)$.
We are particularly interested in $F_{g}(\bold m)$ when $|\bold
m|=g-2$, since they determine relations in $R^{g-2}(\mathcal M_g)$
\begin{equation}\label{eqclass}
\kappa(\bold m)=\frac{(2g-3)!!\bold m! F_{g}(\bold
m)}{2g-2}\kappa_{g-2}.
\end{equation}
It is important to notice that we may extend $F_{g}(\bold m)$ to
be defined for all $\bold m\in N^{\infty}$ using the following Lemma
\ref{const2} and Theorem \ref{const5}.

We define constants $\beta_{\bold L}$ and  $\gamma_{\bold L}$ by
$$\sum_{\bold L+\bold{L'}=\bold b}\frac{(-1)^{||\bold L||}\beta_{\bold L}}{\bold{L'}!(2|\bold{L'}|+1)!!}=0,\qquad \bold b\neq0,$$
$$\sum_{\bold L+\bold{L'}=\bold b}\frac{(-1)^{||\bold L||}\gamma_{\bold L}}{\bold{L'}!(2|\bold{L'}|-1)!!}=0,\qquad \bold b\neq0,$$
with the initial values $\beta_{\bold 0}=\gamma_{\bold 0}=1$. We
also denote their reciprocals by
$$\gamma^{-1}_{\bold L}:=\frac{(-1)^{||\bold L||}}{\bold{L}!(2|\bold{L}|-1)!!},\qquad \beta^{-1}_{\bold L}:=\frac{(-1)^{||\bold L||}}{\bold{L}!(2|\bold{L}|+1)!!}.$$

\begin{lemma}\label{const1} Let $n\geq0$ and $|\bold m|\leq g-2$. Then
\begin{equation} \label{eqconst1}
(2g+n-1)F_{g,n+2}(\bold m)=\sum_{\bold L+\bold{L'}=\bold
m}(2g+n-1-2|\bold{L'}|)\gamma_{\bold L} F_{g,n+1}(\bold{L'}).
\end{equation}

\end{lemma}
\begin{proof} We recall the following recursive
formula from Theorem 5.3 of \cite{LX1}
\begin{multline*}
\langle\tau_{d}\tau_{d_0}\prod_{j=1}^n\tau_{d_j}\kappa(\bold
m)\mid\lambda_g\lambda_{g-1}\rangle_g\\
=\sum_{\bold L+\bold{L'}=\bold m}\frac{\gamma_{\bold L}{\bold
m}!}{\bold{L'}!}\frac{(2d+2d_0+2|\bold L|-1)!!}{(2d-1)!!(2d_0-1)!!}
\langle\tau_{d_0+d+|\bold L|-1}\prod_{j=1}^n\tau_{d_j}\kappa(\bold{L'})\mid\lambda_g\lambda_{g-1}\rangle_g\\
+\sum_{\bold L+\bold{L'}=\bold m}\sum_{j=1}^n\frac{\gamma_{\bold
L}{\bold m}!}{\bold{L'}!}\frac{(2d+2d_j+2|\bold
L|-3)!!}{(2d-1)!!(2d_j-3)!!}\langle\tau_{d_0}\tau_{d_j+d+|\bold
L|-1}\prod_{i\neq
j}\tau_{d_i}\kappa(\bold{L'})\mid\lambda_g\lambda_{g-1}\rangle_g.
\end{multline*}
When $\bold m=0$, it was obtained by Getzler and Pandharipande
\cite{GP} from degree $0$ Virasoro constraints for $\mathbb P^2$.

So we have
\begin{align*}
(2g+n-1)F_{g,n+2}(\bold m)=&\sum_{\bold L+\bold{L'}=\bold
m}(2d+2d_0+2|\bold
L|-1)\gamma_{\bold L}F_{g,n+1}(\bold{L'})\\
&+\sum_{\bold L+\bold{L'}=\bold m}\sum_{j=1}^n(2d_j-1)\gamma_{\bold
L}F_{g,n+1}(\bold{L'})\\
=&\sum_{\bold L+\bold{L'}=\bold m}(2g+n-1-2|\bold{L'}|)\gamma_{\bold
L} F_{g,n+1}(\bold{L'}).
\end{align*}
In the last equation, we used $d+d_0+\sum_{j=1}^n d_j=g+n-|\bold m|$.
\end{proof}
Lemma \ref{const1} can also be proved using the following identity instead
$$\sum_{\bold L+\bold L'=\bold
b}(-1)^{||\bold L||}\binom{\bold b}{\bold L}\langle\tau_{|\bold
L|}\prod_{j=1}^n\tau_{d_j}\kappa(\bold
L')\mid\lambda_g\lambda_{g-1}\rangle_g=\sum_{j=1}^n\langle\tau_{d_j-1}\prod_{i\neq
j}\tau_{d_i}\kappa(\bold b)\mid\lambda_g\lambda_{g-1}\rangle_g,$$
which can be proved by the same argument of Proposition 3.1 of \cite{LX1}.

\begin{lemma}\label{const2} Let $n\geq0$ and $|\bold m|\leq g-2$. Then
\begin{equation} \label{eqconst2}
F_{g,n}(\bold m)=\sum_{\bold L+\bold{L'}=\bold m}\beta^{-1}_{\bold
L} F_{g,n+1}(\bold{L'}).
\end{equation}
\end{lemma}

\begin{proof} We use the following identity
$$(2g-2+n)\langle\prod_{j=1}^n\tau_{d_j}\kappa(\bold
m)\mid\lambda_g\lambda_{g-1}\rangle_g=\sum_{\bold L+\bold L'=\bold
m}(-1)^{||\bold L||}\binom{\bold m}{\bold L}\langle\tau_{|\bold
L|+1}\prod_{j=1}^n\tau_{d_j}\kappa(\bold
L')\mid\lambda_g\lambda_{g-1}\rangle_g,$$ which can be proved by the same argument of Proposition 3.1 of
\cite{LX1}. By a direct calculation as Lemma \ref{const1}, we get the desired result.
\end{proof}

\begin{lemma}\label{const3} Let $n\geq0$ and $|\bold m|\leq g-2$. Then
\begin{equation} \label{eqconst3}
2|\bold m| F_{g,n+1}(\bold m)=\sum_{\substack{\bold e+\bold f+\bold
L =\bold m\\\bold L\neq \bold m}}\beta^{-1}_{\bold e}\gamma_{\bold
f}(2g+n-1-2|\bold L|) F_{g,n+1}(\bold L).
\end{equation}
\end{lemma}
\begin{proof}
The result follows by applying Lemma \ref{const1} to the right hand
side of equation \eqref{eqconst2} in Lemma \ref{const2}. We have
\begin{multline*}
(2g+n-1) F_{g,n+1}(\bold m)=\sum_{\bold e+\bold{b}=\bold m}\beta^{-1}_{\bold
e} (2g+n-1) F_{g,n+2}(\bold{b})\\
=\sum_{\bold e+\bold{b}=\bold m}\beta^{-1}_{\bold
e} \sum_{\bold f+\bold
L =\bold b}\gamma_{\bold
f}(2g+n-1-2|\bold L|) F_{g,n+1}(\bold L)\\
=\sum_{\bold e+\bold f+\bold
L =\bold m}\beta^{-1}_{\bold e}\gamma_{\bold
f}(2g+n-1-2|\bold L|) F_{g,n+1}(\bold L).
\end{multline*}
So we get the desired identity.
\end{proof}

\begin{proposition}\label{const4} Let $n\geq0$ and $|\bold m|\leq g-2$. Then
\begin{equation} \label{eqconst4}
2|\bold m|F_{g,n}(\bold m)=(2g+n-2)\sum_{\substack{\bold L+\bold{L'}=\bold m\\
\bold{L}\neq\bold 0}}C_{\bold L}F_{g,n}(\bold{L'}),
\end{equation}
where $C_{\bold 0}=-1$ and
$$C_{\bold
L}=\sum_{\bold e+\bold f=\bold L}2|\bold e|\beta_{\bold
e}\beta^{-1}_{\bold f}=-\sum_{\bold e+\bold f=\bold
L}\gamma^{-1}_{\bold e}\beta_{\bold f},\qquad \bold L\neq0.$$
\end{proposition}

\begin{proof}
From Lemma \ref{const2} and Lemma \ref{const3}, we have

\begin{multline*}
2|\bold m|\sum_{\bold L+\bold{L'}=\bold m}\beta_{\bold L}
F_{g,n}(\bold{L'})=2|\bold m|F_{g,n+1}(\bold
m)\\=\sum_{\substack{\bold e+\bold f+\bold L =\bold m\\\bold L\neq
\bold m}}\beta^{-1}_{\bold e}\gamma_{\bold f}(2g+n-1-2|\bold
L|)\sum_{\bold {L'}+\bold{L''}=\bold L}\beta_{\bold
{L'}}F_{g,n}(\bold {L''})\\=\sum_{\bold e+\bold f+\bold L =\bold
m}\beta^{-1}_{\bold e}\gamma_{\bold f}(2g+n-1-2|\bold L|)\sum_{\bold
{L'}+\bold{L''}=\bold L}\beta_{\bold {L'}}F_{g,n}(\bold
{L''})-(2g+n-1-2|\bold m|)\sum_{\bold L+\bold{L'}=\bold
m}\beta_{\bold L}F_{g,n}(\bold {L'}).
\end{multline*}

Subtract $2|\bold m|\sum_{\bold L+\bold{L'}=\bold m}\beta_{\bold L}
F_{g,n}(\bold{L'})$ from each side, we have

\begin{multline}\label{eqconst6}
0=(2g+n-1)\sum_{\bold f+{\bold L''}=\bold m}\gamma_{\bold
f}F_{g,n}(\bold {L''})-2\sum_{\bold e+\bold
f+\bold{L'}+\bold{L''}=\bold m}|\bold{L''}|\beta^{-1}_{\bold
e}\gamma_{\bold f}\beta_{\bold {L'}}F_{g,n}(\bold {L''})
\\-2\sum_{\bold e+\bold f+\bold{L'}+\bold{L''}=\bold
m}|\bold{L'}|\beta^{-1}_\bold{e}\gamma_{\bold
f}\beta_{\bold{L'}}F_{g,n}(\bold{L''})-(2g+n-1)\sum_{\bold L+\bold
{L'}=\bold m}\beta_{\bold L}F_{g,n}(\bold{L'}).
\end{multline}

Now we simplify the third term in \eqref{eqconst6}
\begin{multline*}
-\sum_{\bold e+\bold f+\bold{L'}+\bold{L''}=\bold
m}2|\bold{L'}|\beta^{-1}_\bold{e}\gamma_{\bold
f}\beta_{\bold{L'}}F_{g,n}(\bold{L''})=\sum_{\bold e+\bold
f+\bold{L'}+\bold{L''}=\bold
m}2|\bold{e}|\beta_\bold{e}\gamma_{\bold
f}\beta_{\bold{L'}}F_{g,n}(\bold{L''})\\=\sum_{\bold e+\bold
f+\bold{L'}+\bold{L''}=\bold
m}(2|\bold{e}|+1)\beta^{-1}_\bold{e}\gamma_{\bold
f}\beta_{\bold{L'}}F_{g,n}(\bold{L''})-\sum_{\bold e+\bold
f+\bold{L'}+\bold{L''}=\bold m}\beta^{-1}_\bold{e}\gamma_{\bold
f}\beta_{\bold{L'}}F_{g,n}(\bold{L''})\\=\sum_{\bold e+\bold
f+\bold{L'}+\bold{L''}=\bold m}\gamma^{-1}_\bold{e}\gamma_{\bold
f}\beta_{\bold{L'}}F_{g,n}(\bold{L''})-\sum_{\bold L+\bold
{L'}=\bold m}\gamma_{\bold L}F_{g,n}(\bold {L'})\\=\sum_{\bold
L+\bold {L'}=\bold m}\beta_{\bold L}F_{g,n}(\bold {L'})-\sum_{\bold
L+\bold {L'}=\bold m}\gamma_{\bold L}F_{g,n}(\bold {L'}).
\end{multline*}

Substitute into \eqref{eqconst6}, we have

\begin{multline*}
0=(2g+n-1)\sum_{\bold L+\bold {L'}=\bold m}\gamma_{\bold
L}F_{g,n}(\bold {L'})-2\sum_{\bold L+\bold {L'}=\bold
m}|\bold{L'}|\gamma_{\bold L}F_{g,n}(\bold{L''})\\+\sum_{\bold
L+\bold {L'}=\bold m}\beta_{\bold L}F_{g,n}(\bold {L'})-\sum_{\bold
L+\bold {L'}=\bold m}\gamma_{\bold L}F_{g,n}(\bold
{L'})-(2g+n-1)\sum_{\bold L+\bold {L'}=\bold m}\beta_{\bold
L}F_{g,n}(\bold {L'}).
\end{multline*}

So we get

\begin{multline*}
2\sum_{\bold L+\bold {L'}=\bold m}|\bold{L'}|\gamma_{\bold
L}F_{g,n}({L'})=(2g+n-2)\sum_{\bold L+\bold {L'}=\bold
m}\gamma_{\bold L}F_{g,n}(\bold {L'})-(2g+n-2)\sum_{\bold
L+\bold{L'}=\bold m}\beta_{\bold L}F_{g,n}(\bold{L'}).
\end{multline*}

Convoluting both sides by $\gamma^{-1}_{\bold L}$, we have

\begin{multline*}
2|\bold m|F_{g,n}(\bold m)=(2g+n-2)F_{g,n}(\bold
m)-(2g+n-2)\sum_{\bold L+\bold {L'}=\bold m}(\sum_{\bold e+\bold
f=\bold L}\gamma^{-1}_{\bold e}\beta_{\bold f})F_{g,n}(\bold
{L'})\\=-(2g+n-2)\sum_{\substack{\bold L+\bold{L'}=\bold m\\
\bold{L}\neq\bold 0}}(\sum_{\bold e+\bold
f=\bold L}\gamma^{-1}_{\bold e}\beta_{\bold f})F_{g,n}(\bold {L'})\\=-(2g+n-2)\sum_{\substack{\bold L+\bold{L'}=\bold m\\
\bold{L}\neq\bold 0}}(\sum_{\bold e+\bold
f=\bold L}(2|\bold e|+1)\beta^{-1}_{\bold e}\beta_{\bold f})F_{g,n}(\bold {L'})\\=(2g+n-2)\sum_{\substack{\bold L+\bold{L'}=\bold m\\
\bold{L}\neq\bold 0}}(\sum_{\bold e+\bold f=\bold L}2|\bold
e|\beta_{\bold e}\beta^{-1}_{\bold f})F_{g,n}(\bold {L'}).
\end{multline*}

\end{proof}

Our first main result follows as an immediately corollary of the above proposition.

\begin{theorem}\label{const5} Let $|\bold m|\leq g-2$. Then
\begin{equation} \label{eqconst5}
|\bold m|F_g(\bold m)=(g-1)\sum_{\substack{\bold L+\bold{L'}=\bold m\\
\bold{L}\neq\bold 0}}C_{\bold L}F_g(\bold{L'}).
\end{equation}
The constant $C_{\bold L}$ is defined in Proposition \ref{const4}. In particular,
$$F_g(\bm{\delta}_k)=\frac{2g-2}{(2k+1)!!},\qquad k\geq1.$$
$\bm{\delta}_k$ denotes the sequence with $1$ at the $k$-th place
and zeros elsewhere.
\end{theorem}
\begin{proof}
For the last assertion, just note that
$$C_{\bm{\delta}_k}=2k\cdot\beta_{\bm{\delta}_k}=\frac{2k}{(2k+1)!!}.$$
\end{proof}

Theorem \ref{const5} gives an efficient way to compute constants
$F_g(\bold m)$ when $|\bold m|=g-2$, hence tautological relations in $\mathcal
R^{g-2}(\mathcal M_g)$ through equation \eqref{eqclass}.
The identity \eqref{eqconst5} may be expanded to get a closed
formula of $F_g(\bold m)$.

\begin{corollary}\label{const6} For $\bold m\neq\bold 0$, we have
\begin{equation*}F_g(\bold m)
=\sum_{k=1}^{||\bold m||}(g-1)^k\sum_{\substack {\bold
m=\bold{m_1}+\cdots+\bold{m_k}\\\bold{m_i}\neq\bold
0}}\frac{\prod_{j=1}^k
C_{\bold{m_j}}}{\prod_{j=1}^k|\bold{m_1}+\cdots+\bold{m_j}|}.
\end{equation*}
\end{corollary}

Now we come to our second main result, in contrast with Theorem \ref{const5},  it is a formula that gives
recursive relations only among those $F_g(\bold m)$ with $|\bold m|=g-2$.
\begin{theorem} \label{recur1}
Let $g>2$ and $|\bold m|=g-2$. Then we have
\begin{equation} \label{eqconst8}
(||\bold m||-1)F_{g}(\bold m)=\sum_{\substack{\bold L+\bold{L'}=\bold m\\||\bold{L'}||\geq2}}D_{g,\bold{L'}}
\frac{(\bold L+\bm\delta_{g-2-|\bold{L}|})!}{\bold L!}F_{g}(\bold L+\bm\delta_{g-2-|\bold{L}|}),
\end{equation}
where the constant $D_{g,\bold{L'}}$ is given by
$$D_{g,\bold{L'}}=-\frac{1}{\bold{L'}!}+
\frac{2g-1}{2(g-2)}\sum_{\substack{\bold{L_1}+\bold{L_2}=\bold{L'} \\||\bold{L_1}||\geq1}}C_{\bold{L_1}}\frac{(1+2|\bold{L_1}|)!!}{\bold{L_2}!}.$$
The constant $C_{\bold L}$ is defined in Proposition \ref{const4}.
\end{theorem}
\begin{proof}
By the projection formula and $\kappa_0=2g-2$, we have
\begin{multline}
\langle\tau_1\kappa(\bold m)\mid \lambda_{g}\lambda_{g-1}\rangle_g=\int_{\overline\sM_g}\pi_*(\psi\cdot\prod_{i\geq1}(\pi^*\kappa_i+\psi^i)^{m_i})
 \lambda_{g}\lambda_{g-1}\\
=\sum_{\bold L+\bold{L'}=\bold m}\binom{\bold m}{\bold L}\int_{\overline\sM_g}\kappa(\bold L)\kappa_{|\bold{L'}|} \lambda_{g}\lambda_{g-1}\\
=(2g-2)\int_{\overline\sM_g}\kappa(\bold m) \lambda_{g}\lambda_{g-1}+
\sum_{\substack{\bold L+\bold{L'}=\bold m\\\bold{L'}\ne\bold0}}\binom{\bold m}{\bold L}\int_{\overline\sM_g}\kappa(\bold L)\kappa_{|\bold{L'}|} \lambda_{g}\lambda_{g-1}.
\end{multline}

From equation \eqref{eqconst7}, the above equation becomes
\begin{multline}\label{eqconst11}
F_{g,1}(\bold m)=F_g(\bold m)+\frac{1}{2g-2}\sum_{\substack{\bold
L+\bold{L'}=\bold m\\\bold{L'}\ne\bold0}}\frac{1}{\bold
L!\bold{L'}!}
F_{g}(\bold L+\bm\delta_{|\bold{L'}|})(\bold L+\bm\delta_{|\bold{L'}|})!\\
=(1+\frac{||\bold m||}{2g-2})F_g(\bold m)+\frac{1}{2g-2}\sum_{\substack{\bold L+\bold{L'}=\bold m\\||\bold{L'}||\geq2}}\frac{1}{\bold L!\bold{L'}!}
F_{g}(\bold L+\bm\delta_{|\bold{L'}|})(\bold L+\bm\delta_{|\bold{L'}|})!
\end{multline}

Take $n=1$ in Proposition \ref{const4}, we have
$$2|\bold m|F_{g,1}(\bold m)=(2g-1)\sum_{\substack{\bold L+\bold{L'}=\bold m\\
\bold{L}\neq\bold 0}}C_{\bold L}F_{g,1}(\bold{L'})$$ Substitute
equation \eqref{eqconst11} into both sides of the above identity and
adjust the indices, we get
\begin{multline*}
(g-2)(||\bold m||-1)F_{g}(\bold m)\\=\frac{2g-1}{2}\sum_{\substack{\bold L+\bold{L'}=\bold m\\||\bold{L'}||\geq2}}
\sum_{\substack{\bold{L_1}+\bold{L_2}=\bold{L'}\\||\bold{L_1}||\geq\bold1}}C_{\bold{L_1}}\frac{(1+2|\bold{L_1}|)!!}{\bold L!\bold{L_2}!}
F_{g}(\bold L+\bm\delta_{g-2-|\bold{L}|})(\bold L+\bm\delta_{g-2-|\bold{L}|})!\\
-(g-2)\sum_{\substack{\bold L+\bold{L'}=\bold m\\||\bold{L'}||\geq2}}\frac{1}{\bold L!\bold{L'}!}
F_{g}(\bold L+\bm\delta_{g-2-|\bold{L}|})(\bold L+\bm\delta_{g-2-|\bold{L}|})!,
\end{multline*}
which is just our desired identity.
\end{proof}

We know $F_g(\bm\delta_{g-2})=\frac{2g-2}{(2g-3)!!}$. So
Theorem \ref{recur1} gives a recursive formula to compute $F(\bold m)$ when $|\bold m|=g-2$
by induction on $||\bold m||$. We have written a maple program to implement Theorems \ref{const5}
and \ref{recur1}, both of them give the correct results.

\begin{remark}
Here we comment on the signs of the constants that we met in this
section. Although we are not able to give a proof for now, numerical
evidence strongly suggests that $\beta_{\bold L} >0, \gamma_{\bold
L}>0$ for all $L \in N^\infty$ and $C_{\bold L}>0, D_{g,\bold L}>0$
for all. $\bold L \neq 0$ and $g\geq 3$. When $\bold L=(1^k)$, this
is easy to verify (see the next section). One may wonder at the
seemly exceptional initial values $C_0=D_0=-1$; in fact, their
negativity is essential for the tautological relation in Theorem
\ref{tauto}.
\end{remark}

\vskip 30pt
\section{Identities of Bernoulli numbers} \label{Bern}

Let $k\geq0$. Denote by $C_k$ the constants $C_{(k,0,0,\dots)}$
defined in Proposition \ref{eqconst4}. The same convention apply for
$\beta_k,\beta^{-1}_k,\gamma_k,\gamma^{-1}_k, D_{g,k}$. Recall that
we have (see \cite{LX1})
\begin{align}
\beta_{k}&=\frac{(-1)^{k}(2-2^{2k})B_{2k}}{k!(2k-1)!!},\label{eqbeta1}\\
\beta^{-1}_k&=\frac{(-1)^k}{k!(2k+1)!!},\label{eqbeta2}\\
\gamma_{k}&=\frac{E_{2k}}{(2k-1)!!},\label{eqgamma1}\\
\gamma^{-1}_{k}&=\frac{(-1)^k}{k!(2k-1)!!},\label{eqgamma2}
\end{align}
where $B_k, E_k$ are the Bernoulli and Euler numbers respectively
$$\frac{t}{e^t-1}=\sum_{m=0}^{\infty}B_m\frac{t^m}{m!},$$
$$\sec t=\sum_{m=0}^\infty E_m\frac{t^{m}}{m!}.$$
We have $B_1=-\frac12$, $B_{2k+1}=0$ for $k\geq1$ and $E_{2k+1}=0$
for $k\geq0$.

We first prove a simple closed formula for $C_k$.
\begin{lemma}\label{const7} Let $k\geq0$. Then
$$C_k=-\sum_{j=0}^k\gamma^{-1}_j\beta_{k-j}=\frac{(-1)^{k+1} 2^{2k} B_{2k}}{k!(2k-1)!!}.$$
\end{lemma}
\begin{proof}
From \eqref{eqbeta1} and \eqref{eqgamma2}, we have
\begin{multline*}
-\sum_{j=0}^k\gamma^{-1}_j\beta_{k-j}=\sum_{j=0}^k\frac{(-1)^{k-j}}{(k-j)!(2k-2j-1)!!}\cdot\frac{(-1)^{j+1}(2-2^{2j})B_{2j}}{j!(2j-1)!!}\\
=\frac{(-1)^{k+1}(2-2^{2k})B_{2k}}{k!(2k-1)!!}+\sum_{j=0}^{k-1}\frac{(-1)^{k-j}}{(k-j)!(2k-2j-1)!!}\cdot\frac{(-1)^{j+1}(2-2^{2j})B_{2j}}{j!(2j-1)!!}.
\end{multline*}
So we need to prove
$$2(2^{2k}-1)(2k+1)B_{2k}=\sum^{k-1}_{j=0}\binom{2k+1}{2j}(2-2^{2j})(2k-2j+1)B_{2j},$$
which is easily seen to follow from the following lemma.
\end{proof}

\begin{lemma}Let $n\geq2$. Then
$$2(2^n-1)(n+1)B_n=\sum^{n-1}_{j=0}\binom{n+1}{j}(2-2^j)(n-j+1)B_j.$$
\end{lemma}
\begin{proof} Multiplying $\frac{1}{(n+1)!}$ on both sides and taking generating functions, we have
\begin{align*}
A:=\sum_{n=0}^{\infty}\frac{2(2^n-1)B_n}{n!}t^n&=\frac{-2t}{\sinh t}+2\sum_{n=0}^{\infty}\frac{B_n}{n!}t^n\\
&=\frac{-2t}{\sinh t}+t\coth(t/2)-t
\end{align*}
and

\begin{align*}
B:=&\sum_{n=0}^{\infty}\sum^{n-1}_{j=0}\frac{2-2^j}{j!(n+1-j)!}(n-j+1)B_j t^n\\
=&\sum_{n=0}^{\infty}\sum^{n}_{j=0}\frac{2-2^j}{j!(n-j)!}B_j
t^n-\sum_{n=0}^{\infty}\frac{2-2^n}{n!}B_nt^n\\
=&\frac{t e^t}{\sinh t}-\frac{t}{\sinh t}.
\end{align*}

Since the lemma only cares about coefficients of $t^n$ with
$n\geq2$, it is sufficient to prove that $$A+2t=B.$$
Multiply
$(e^{2t}-1)/t$ on each side,
it is easy to check that both sides
equal $2e^{2t}-2e^t$.
\end{proof}

We denote $F_g(k,0,0,\dots)$ by $F_g(k)$.

\begin{lemma} \label{Fg} Let $k\geq1$. Then
\begin{equation}\label{eqFg}
F_g(k)=2^{2k}\sum_{n=1}^{k}(g-1)^n\sum_{\substack{a_1+\cdots+a_n=k\\a_i>0}}
\prod_{j=1}^n\frac{|B_{2a_j}|}{a_j!(2a_j-1)!!\cdot|a_1+\cdots+a_j|}.
\end{equation}
\end{lemma}
\begin{proof} This follows directly from Corollary \ref{const6} and Lemma
\ref{const7}. Also note that $|B_{2m}|=(-1)^{m+1} B_{2m},\ m>0$.
\end{proof}

\begin{lemma}\label{Fzero} Let $k\geq0$. Then
$$F_0(k)=\frac{(-1)^k 2^{2k}}{(k+1)!(2k+1)!!}.$$
\end{lemma}

\begin{proof}
By Lemma \ref{const7} and Theorem \ref{const5}, we need to prove
that for $k\geq0$,
\begin{align*}
\frac{(-1)^kk\cdot2^{2k}}{(k+1)!(2k+1)!!}=&\sum_{j=1}^k\frac{(-1)^{j+1}2^{2j}B_{2j}}{(2j-1)!!j!}\cdot
\frac{(-1)^{k-j}2^{2(k-j)}}{(k-j+1)!(2k-2j+1)!!}\\
=&\frac{(-1)^{k+1}2^{2k}}{(k+1)!(2k+1)!!}+\sum_{j=0}^k\frac{(-1)^{j+1}2^{2j}B_{2j}}{(2j-1)!!j!}\cdot
\frac{(-1)^{k-j}2^{2(k-j)}}{(k-j+1)!(2k-2j+1)!!}
\end{align*}

It is not difficult to simplify the above equation to
$$\sum_{j=0}^k\binom{2k+2}{2j}B_{2j}=k+1.$$

We have
\begin{align*}
\sum_{j=0}^k\binom{2k+2}{2j}B_{2j}=&\sum_{j=0}^{2k+1}\binom{2k+2}{j}B_j-\binom{2k+2}{1}B_1\\
=&k+1,
\end{align*}
where we used the well-known formula
$$\sum_{j=0}^m\binom{m+1}{j}B_j=0,\qquad m>0.$$
This completes the proof.
\end{proof}

The following Faber-Zagier's formula \cite{Fa} may be proved using
the Faber intersection number conjecture and the Cauchy residue
formula \cite{Zh}.

\begin{proposition}{\bf (Faber-Zagier)} Let $g\geq2$. Then
\begin{equation}\label{eqFaZa}
F_g(g-2)=\frac{2^{2g-4}(g-2)!}{(2g-3)!!}.
\end{equation}
\end{proposition}
\begin{proof}
By \eqref{eqkmz} and \eqref{eqclass}, it is easy to see that
Faber-Zagier's formula is equivalent to the following combinatorial
lemma.
\end{proof}

The proof of the following Lemma \ref{zhou} is due to Jian Zhou
\cite{Zh}.
\begin{lemma} \label{zhou}
Let $g\geq 1$. Then
$$\sum_{k=1}^g\left(\frac{(-1)^k}{k!}(2g+1+k)\sum_{g=m_1+\dots+m_k\atop m_i>0}\binom{2g+k}{2m_1+1,\dots,2m_k+1}\right)=(-1)^g2^{2g}(g!)^2$$

\end{lemma}
\begin{proof}
We will use Cauchy's residue formula
\begin{align*}
&\sum^g_{k=1}\frac{(-1)^k}{k!}(2g+1+k)\sum_{m_1+\dots+m_k=g\atop
m_i>0}\binom{2g+k}{2m_1+1,\dots,2m_k+1}\\
=&\sum_{k=1}^g(-1)^k\prod_{j-1}^{2g+1}(k+j)\sum_{m_1+\dots+m_k=g\atop
m_i>0}\frac{1}{(2m_1+1)!}\dots\frac{1}{(2m_k+1)!}\\
=&\sum_{k=1}^g(-1)^k\prod_{j-1}^{2g+1}(k+j)\Res_{z=0}\frac{1}{z^{2g+1+k}}f(z)^k\\
=&\Res_{z=0}\left(\sum_{k=1}^g(-1)^k\prod_{j-1}^{2g+1}(k+j)\frac{1}{z^{2g+1+k}}f(z)^k\right),
\end{align*}
where
$$f(z)=\sum_{m>0}\frac{z^{2m+1}}{(2m+1)!}.$$

To take the summation over $k$, we notice that
$$\prod_{j-1}^{2g+1}(k+j)\frac{1}{z^{2g+1+k}}=-z\partial_w^{2g+1}|_{w=z}w^{-k-1},$$
hence we can proceed as follows:
\begin{align*}
&\Res_{z=0}\sum_{k=1}^g(-1)^k\prod_{j-1}^{2g+1}(k+j)\frac{1}{z^{2g+1+k}}f(z)^k\\
=&-\Res_{z=0}\sum_{k=1}^g(-1)^kz\partial_w^{2g+1}|_{w=z}w^{-k-1}f(z)^k\\
=&-\Res_{z=0}z\partial_w^{2g+1}|_{w=z}\frac{-w^{-2}f(z)}{1+w^{-1}f(z)}\\
=&-\Res_{z=0}z\partial_w^{2g+1}|_{w=z}(\frac{1}{w+f(z)}-\frac{1}{w})\\
=&(2g+1)!\Res_{z=0}\frac{z}{(z+f(z))^{2g+2}}\\
=&(2g+1)!\Res_{z=0}\frac{z}{\sinh^{2g+2}z}\\
=&(2g+1)!\frac{1}{2\pi i}\oint_{z}\frac{z}{\sinh^{2g+2}z}dz.
\end{align*}

To evaluate the contour integral, we make the following change of
variable. Take $u=\sinh z$, $z=\arcsinh u$, so that
$dz=(1+u^2)^{-1/2}du$, therefore,$$\frac{1}{2\pi
i}\oint_{z}\frac{z}{\sinh^{2g+2}z}dz=\frac{1}{2\pi
i}\oint_{u}\frac{1}{u^{2g+2}}\arcsinh u\cdot(1+u^2)^{-1/2}du$$

However, if we write $$\arcsinh
u\cdot(1+u^2)^{-1/2}=\sum^{\infty}_{g=0}a_gu^{2g+1},$$ and
differentiate both sides
\begin{align*}
\sum^{\infty}_{g=0}(2g+1)a_gu^{2g}&=-u(1+u^2)^{-3/2}\arcsinh
u+(1+u^2)^{-1}\\
&=(1+u^2)^{-1}(1-\sum^{\infty}_{g=0}a_gu^{2g+2}),
\end{align*}

We get $a_0=1$ and
\begin{equation*}
a_{g+1}=-\frac{2g+2}{2g+3}a_g, \ g\geq0.
\end{equation*}

From the recursion relation one gets the following unique solution:
$$a_{g}=-\frac{(-1)^g}{(2g+1)!}2^{2g}(g!)^2.$$

This completes the proof.
\end{proof}

By Lemma \ref{Fg}, we may rewrite Faber-Zagier's formula
\eqref{eqFaZa} as the following identity of Bernoulli numbers.
\begin{proposition} Let $g\geq1$. Then
\begin{equation}\label{eqbern}
\sum_{n=1}^{g}(g+1)^n\sum_{\substack{a_1+\cdots+a_n=g\\a_i>0}}
\prod_{j=1}^n\frac{|B_{2a_j}|}{a_j!(2a_j-1)!!\cdot|a_1+\cdots+a_j|}=\frac{g!}{(2g+1)!!}.
\end{equation}
\end{proposition}

It shall be interesting to find a direct combinatorial proof of
\eqref{eqbern}. However, we are not able to find a general explicit
expression for $F_g(k)$ and it seems not easy to prove
\eqref{eqbern} using Theorem \ref{const5}.

For comparison, we note that Lemma \ref{Fg} and Lemma \ref{Fzero}
give the following identity
\begin{equation}\label{eqbern2}
\sum_{n=1}^{k}\sum_{\substack{a_1+\cdots+a_n=k\\a_i>0}}
\prod_{j=1}^n\frac{B_{2a_j}}{a_j!(2a_j-1)!!\cdot|a_1+\cdots+a_j|}=\frac{1}{(k+1)!(2k+1)!!}.
\end{equation}

\vskip 30pt
\section{Proof of Theorem \ref{tauto}}
First note that the explicit value of $D_{g,k}$ in \eqref{eqconst10}
follows from Lemma \ref{const7}. Now Theorem \ref{tauto} is an easy
consequence of Theorem \ref{recur1} and Faber-Zagier's formula
\eqref{eqfaza2},
\begin{align*}
\sum_{i=0}^{g-4}D_{g,g-2-i}\frac{\kappa_{1}^i\kappa_{g-2-i}}{i!}
=&\frac{(2g-3)!!}{2g-2}\sum_{i=0}^{g-4}D_{g,g-2-i}F_g(1^i, g-2-i)\\
=&\frac{(2g-3)!!}{2g-2}(g-3)F_g(1^{g-2})\\
=&\frac{g-3}{(g-2)!}\kappa_{1}^{g-2}.
\end{align*}
On the other hand, from $D_{g,0}=-1,\ D_{g,1}=\frac{g+1}{g-2}$, we
have
\begin{align*}
D_{g,0}\frac{\kappa_{1}^{g-2}\kappa_0}{(g-2)!}+D_{g,1}\frac{\kappa_{1}^{g-2}}{(g-3)!}
=&-\frac{2g-2}{(g-2)!}\kappa_{1}^{g-2}+\frac{g+1}{(g-2)!}\kappa_{1}^{g-2}\\
=&\frac{3-g}{g-1}\kappa_1^{g-2}.
\end{align*}
Adding up the above two sets of equations, we complete the proof of
Theorem \ref{tauto}.

\begin{proposition} Let $g\geq3$. If we assume that Faber's perfect pairing
conjecture is true in codimension one, namely the natural product
$$R^1(\mathcal M_g) \times R^{g-3}(\mathcal M_g)
\rightarrow R^{g-2}(\mathcal M_g) \cong \mathbb Q$$ is
nondegenerate, then the following relation
\begin{equation}\label{eqconst9}
D_{g,
g-2}\frac{g-1}{2^{2g-5}((g-2)!)^2}\kappa_1^{g-3}+\sum_{i=0}^{g-3}D_{g,g-3-i}\frac{\kappa_1^{i}
\kappa_{g-3-i}}{(i+1)!}=0
\end{equation}
holds in $R^{g-3}(\mathcal M_g)$.
\end{proposition}
\begin{proof}
By Theorem \ref{tauto} and Faber-Zagier's formula \eqref{eqfaza2},
we have
\begin{multline*}
\sum_{i=0}^{g-2}D_{g,g-2-i}\frac{\kappa_1^i
\kappa_{g-2-i}}{i!}=D_{g,g-2}\kappa_{g-2}+\sum_{i=1}^{g-2}D_{g,g-2-i}\frac{\kappa_1^i
\kappa_{g-2-i}}{i!}\\=\kappa_1\cdot\left(D_{g,
g-2}\frac{g-1}{2^{2g-5}((g-2)!)^2}\kappa_1^{g-3}+\sum_{i=1}^{g-2}D_{g,g-2-i}\frac{\kappa_1^{i-1}
\kappa_{g-2-i}}{i!}\right)=0.
\end{multline*}

So the relation \eqref{eqconst9} follows from Faber's perfect
pairing conjecture in codimension one and a substitution of $i$ by
$i+1$.
\end{proof}

\begin{example}
Take $g=6$ in \eqref{eqconst9}, we get \begin{align*}
&D_{6,4}\frac{5}{2^7(4!)^2}\kappa_1^3+D_{6,0}\frac{10}{4!}\kappa_1^3+D_{6,1}\frac{1}{3!}\kappa_1^3
+D_{6,2}\frac{1}{2!}\kappa_1\kappa_2+D_{6,3}\kappa_3
\\=&-\frac{5161}{41472}\kappa_1^3+\frac{49}{24}\kappa_1\kappa_2+\frac{395}{72}\kappa_3
\end{align*}

Substitute into the above identity the following relations in
$R^*(\mathcal M_6)$ computed by Faber \cite{Fa},
$$\kappa_3=\frac{5}{2304}\kappa_1^3, \qquad \kappa_1\kappa_2=\frac{127}{2304}\kappa_1^3
,$$ we then checked that the tautological relation \eqref{eqconst9}
holds when $g=6$.
\end{example}

\begin{remark}
Counterexamples of analogues of Faber's perfect pairing conjecture
on partially compactified moduli spaces of curves have recently been
found by Cavalieri and Yang \cite{CaYa}, but we still have reasons
to believe that Faber's original perfect pairing conjecture for
$R^*(\mathcal M_g)$ is true, in view of the beautiful combinatorial
structures of $R^*(\mathcal M_g)$ even in the top degree. On the
other hand, finding explicit tautological relations in lower degrees
is a much harder but rewarding problem. In this regard, we refer the
interested readers to \cite{Fa, Io}.
\end{remark}

$$ \ \ \ \ $$

\end{document}